\documentclass{amsart}
\usepackage{fullpage}
\usepackage{amsmath,amsthm,amssymb}
\usepackage{stmaryrd}
\usepackage{tikz}
\usetikzlibrary{matrix}

\newcommand{\Z}{\mathbb{Z}}
\newcommand{\A}{\hat{\mathbb{A}}}
\newcommand{\G}{\hat{\mathbb{G}}}
\newcommand{\C}{\mathbb{C}}
\newcommand{\CP}{\C \mathrm{P}}

\newcommand{\<}{\langle}
\renewcommand{\>}{\rangle}
\newcommand{\from}{\leftarrow}

\newcommand{\sheaf}[1]{\mathcal{#1}}

\DeclareMathOperator{\Alt}{Alt}
\DeclareMathOperator*{\limone}{lim^1}
\DeclareMathOperator{\Spf}{Spf}
\DeclareMathOperator{\Spec}{Spec}
\DeclareMathOperator{\Hom}{Hom}
\DeclareMathOperator{\id}{id}
\DeclareMathOperator{\Tot}{Tot}
\DeclareMathOperator{\Span}{Span}
\DeclareMathOperator{\sch}{sch}
\DeclareMathOperator{\colim}{colim}

\newtheoremstyle{thmstyle}{\topsep}{0.5em}{}{}{\bf}{:}{0.5em}{}
\theoremstyle{thmstyle}
\newtheorem{thm}{Theorem}
\newtheorem{defn}[thm]{Definition}
\newtheorem{cor}[thm]{Corollary}
\newtheorem{lem}[thm]{Lemma}

\begin{document}

\title{The Morava $E$-theory of Eilenberg-Mac Lane spaces}
\author{Eric Peterson}

\maketitle

\begin{abstract}
We deform the Ravenel-Wilson computation of the Morava $K$-homology of Eilenberg-Mac Lane spaces to obtain a similar description of their completed Morava $E$-homology.  This yields both a cohomological description and an interpretation on the level of formal schemes: the scheme associated to the $E$-cohomology of the space $K(\Z/p^\infty, q)$ is the $q$th exterior power of the $p$-divisible group associated to the versal infinitesimally deformed formal group over Lubin-Tate space.
\end{abstract}

\section{Foreword}

One of the most endearing and essential features of topology is the wealth of representability results.  For each machine constructed to study spaces, a space (or space-like object, like a spectrum) can be constructed that in an appropriate sense encodes the original machine, usually in such a way that the representing space inherits many of its algebraic properties.  In another direction, oftentimes when the study of pure algebra dictates the existence of some structure appearing in the results of these topological machines, a great deal of effort goes into lifting versions of these facts to a level intrinsic to topology.  In these ways, machines can be turned on each other and themselves, and after sufficient iterations of the process there is such an overabundance of structure that the answer --- the evaluation of one machine on another --- can't \emph{not} be known.

It is in this context that Ravenel and Wilson~\cite{RW} considered the Morava $K$-theory of Eilenberg-Mac Lane spaces.  Eilenberg-Mac Lane spaces arise from representability in the first sense: they are spaces representing singular cohomology, and as such they come equipped with a graded ring structure.  We also have several good models for them, even interacting well with the ring structure.  Morava $K$-theory arises from representability in the second sense: it is connected to certain geometric points on the moduli of $p$-local formal groups.  As such, its coefficient ring is a graded field; it comes equipped with K\"unneth isomorphisms; and its action on the crucial space $K(\Z, 2)$ is dictated by construction.  This information is precisely the input used by Ravenel and Wilson to make the computation of $K_* K(\Z/p^j, q)$ for all $j$ and $q$.

The deformation theory of formal groups has been studied at length, and the infinitesimal deformations of these geometric points in particular was completed in the 60s by Lubin and Tate~\cite{LT}.  This deformation space is realized in algebraic topology as Morava $E$-theory, and the structure of the $E$-(co)homology of Eilenberg-Mac Lane spaces has been long believed by experts to follow along the same lines as the original result of Ravenel and Wilson --- for instance, it has slipped into the dictionary appendix of Strickland's treatise~\cite{StricklandFPFP} --- but a proof has not appeared.  Providing a proof is the purpose of this note.

\subsection{Outline of proof}

To a ring-valued cohomology theory $R$ and a CW complex $X$, we can associate a formal scheme $X_R := \Spf R^* X = \colim_\alpha \Spec R^* X_\alpha$, where $X_\alpha$ runs over the finite subcomplexes of $X$.  Taking $X$ to be the space $K(\Z/p^\infty, q) := \colim_j K(\Z/p^j, q)$ and $R$ to be a Morava $E$-theory spectrum $E$, we wish to identify $K(\Z/p^\infty, q)_E$ with the formal scheme $\CP^\infty_E[p^\infty]^{\wedge q}$ --- the $q$-fold exterior power of the $p$-divisible group associated to the formal group $\CP^\infty_E$.

The homology $R_* X$ is often more tractable to compute, for various reasons.  When $R^* X$ is free and even-concentrated, homology is related to cohomology in potentially two ways: first, $R_* X$ is always the ring of functions on the internal Hom-scheme $\underline{\Hom}(X_R, \G_m)$, and second, $R_* X$ is sometimes a $R_*$-coalgebra with $\Spf (R_* X)^\vee \cong X_R$, where the system in $\Spf$ runs over the finite dimensional subcoalgebras of $R_* X$.

With these in mind, the proof falls into two steps:
\begin{enumerate}
\item First, we study the coalgebras $E_*^\vee K(\Z/p^\infty, q) := (E^* K(\Z/p^\infty, q))^\vee$, thought of as homology groups.  From a coalgebraic group $M$, we can form its $q$-fold exterior powers $\Alt^q M$, and we show that $E^\vee_* K(\Z/p^\infty, q)$ is isomorphic to $\Alt^q E^\vee_* B\Z/p^\infty$.  Because $E^* B\Z/p^\infty$ is a formal power series ring, the object $\underline{\Hom}((B\Z/p^\infty_E)^{\wedge q}, \G_m)$ makes sense, and this step demonstrates an isomorphism \[\underline{\Hom}((B\Z/p^\infty_E)^{\wedge q}, \G_m) \from \underline{\Hom}(K(\Z/p^\infty, q)_E, \G_m).\]
\item Second, we want to ensure the existence of the scheme $(B\Z/p^\infty_E)^{\wedge q}$ itself and the map $(B\Z/p^\infty_E)^{\wedge q} \to K(\Z/p^\infty, q)_E$ inducing the one above.  Certain coalgebras, those said to have ``good bases'', have an especially good duality with formal schemes, in which case we can show two things: that the coalgebra $E^\vee_* K(\Z/p^\infty, q)$ corresponds to the formal scheme $K(\Z/p^\infty, q)_E$, and that a colimit of coalgebras equates to a colimit of the corresponding formal schemes.  We demonstrate that our coalgebras $\Alt^q E^\vee_* B\Z/p^\infty$ do indeed have this property, and so the isomorphism on homology is the involutive dual of an isomorphism on cohomological formal schemes, described by the same universal properties.  In summary, the morphism \[(B\Z/p^\infty_E)^{\wedge q} \to K(\Z/p^\infty, q)_E\] makes sense and is an isomorphism.
\end{enumerate}

\subsection{Acknowledgements}

The author would like to thank Matthew Ando for suggesting this project and for enduring a long, winding trail of questions, and Constantin Teleman for his hospitality in the summer of 2011.  He acknowledges support from NSF grants DMS-0838434 ``EMSW21MCTP: Research Experience for Graduate Students'' while at UIUC, and later from RTG DMS-0838703 while at UC-Berkeley.

\section{Definitions, notation, and known results}

The symbol $n$ will denote a fixed positive integer throughout, and $p$ will be an odd prime.

\begin{defn} \cite[Section 1.1]{HS}:
$K$ will denote the $n$th Morava $K$-theory spectrum.  Then, define the completed Johnson-Wilson spectrum $E$ as $E = E(n) \wedge_{MU} L_{w_n^{-1} MU} MU$; see~\cite[Section 1.1]{HS}.  This spectrum has homotopy $\pi_* E = \Z[v_1, \ldots, v_{n-1}, v_n^{\pm}]^\wedge_{I}$, where ${}^\wedge_I$ denotes completion at the ideal $I = \< p, v_1, \ldots, v_{n-1} \>$.
\end{defn}

\begin{defn} \cite[Definition 8.3]{HS}:
Define $E^\vee_* X = \pi_* L_{K}(E \wedge X)$.  This is a certain completion of the standard homology functor $E_* X$, but it is itself \emph{not} a homology functor, as $L_K$ is not smashing.
\end{defn}

Hovey and Strickland analyze the properties of $E$ in great detail.  Here are a few crucial results:
\begin{thm} \cite[Proposition 2.5]{HS}:
If $K^* X$ is even-concentrated, then $E^* X$ is even-concentrated, pro-free, and satisfies $E^* X \hat\otimes_{E^*} K^* = K^* X$.
\end{thm}

\begin{thm} \cite[Proposition 8.4]{HS}:
If $K_* X$ is even-concentrated, then $E^\vee_* X$ is even-concentrated, pro-free, and satisfies $E^\vee_* X \hat\otimes_{E_*} K_* = K_* X$.
\end{thm}

\begin{thm} \cite[Propositions 7.10 and 8.4]{HS}:
There is a Milnor exact sequence 
\[ 0 \to \limone_I (E/I)_{m+1} X \to E^\vee_m X \to \lim_I (E/I)_m X \to 0.\]
\end{thm}

Together, these theorems motivate the notation $E^\vee$.  In the situation in which these theorems apply, we have isomorphisms \[E^\vee_* X \cong (E_* X)^\wedge_I \cong (E^* X)^\vee.\]

Now, recall the main results of Ravenel and Wilson.  First, if $E$ and $F$ are ring $\Omega$-spectra and $E$ has enough K\"unneth isomorphisms so that \[E_* (\Omega^{\infty-n} F \times \Omega^{\infty-m} F) \cong E_* \Omega^{\infty-n} F \otimes_{E_*} E_* \Omega^{\infty-m} F,\] then the $E_*$-coalgebras $E_* \Omega^{\infty-*} F$ assemble into a graded coalgebraic ring~\cite[Section 3]{RW}.  We enumerate the various algebraic structures to help fix notation:
\begin{center}
\begin{tabular}{l|l}
Map & Source \\
\hline
$\cdot: E_s \times E_t \Omega^{\infty-n} F \to E_{s+t} \Omega^{\infty-n} F$ & scalar multiplication in $E$-theory, \\
$+: E_s \Omega^{\infty-n} F \times E_s \Omega^{\infty-n} F \to E_s \Omega^{\infty-n} F$ & coaddition on $\mathbb{S}$, \\
$\psi: E_s \Omega^{\infty-n} F \to \bigoplus_t \left( E_t \Omega^{\infty-n} F \otimes_{E_*} E_{s-t} \Omega^{\infty-n} F\right)$ & diagonal on $\Omega^{\infty-n} F$, \\
$\ast: E_s \Omega^{\infty-n} F \otimes_{E_*} E_t \Omega^{\infty-n} F \to E_{s+t} \Omega^{\infty-n} F$ & addition in $F$-theory, \\
$\chi: E_s \Omega^{\infty-n} F \to E_s \Omega^{\infty-n} F$ & negation in $F$-theory, \\
$\circ: E_s \Omega^{\infty-n} F \otimes_{E_*} E_t \Omega^{\infty-m} F \to E_{s+t} \Omega^{\infty-n-m} F$ & multiplication in $F$-theory.
\end{tabular}
\end{center}
These maps satisfy a whole host of properties.  The most complicated such property is distribution of the $\circ$-product over the $\ast$-product: writing $\psi x = \sum x' \otimes x''$, we have \[x \circ (y \ast z) = \sum (x' \circ y) \ast (x'' \circ z).\]  Given a coalgebraic group (i.e., a Hopf algebra), one can construct the free graded coalgebraic ring on it in an obvious way; see Hunton and Turner~\cite{HT} or another document of Ravenel and Wilson~\cite{RW77}.
\begin{thm} \cite[Corollary 11.3]{RW}
The graded coalgebraic ring associated to the Morava $K$-theory of Eilenberg-Mac Lane spaces is free on the coalgebraic group $K_* B\Z/p^j$, i.e., \[\bigoplus_{q=0}^n K_* K(\Z/p^j, q) \cong \bigoplus_{q=0}^\infty \Alt^q K_* B\Z/p^j.\]
\end{thm}

The language of formal geometry has become essential to the study and organization of algebraic topology.  We must at least recall this much:
\begin{defn}
Denote the affine formal $k$-scheme $\Spf k \llbracket x_1, \ldots, x_N \rrbracket$ by $\A^N$.  A formal scheme will be said to be a formal variety (of dimension $N$) if it is isomorphic to $\A^N$.
\end{defn}
\begin{thm} (\cite[pg. 12]{Demazure})
For $k$ a field, there is an equivalence of categories between $k$-coalgebras and formal $k$-schemes.
\end{thm}
\begin{cor}
If $k$ is an even-concentrated field spectrum and $X$ is a space with even-concentrated cohomology, then the coalgebra $k_* X$ and the formal scheme $\Spf k^* X$ are equivalent under the above equivalence of categories.
\end{cor}

In this language, the Ravenel-Wilson result is expressed as:
\begin{cor} \cite[Theorem 12.4]{RW}
Dually, we have an isomorphism of graded ring schemes \[\bigoplus_{q=0}^\infty \Spf K^* K(\Z/p^j, q) \cong \bigoplus_{q=0}^\infty \left(\Spf K^* B\Z/p^j\right)^{\wedge q}.\]  Each graded component is a formal variety.
\end{cor}

\section{$E^\vee$-Homology computation}

The Hovey-Strickland results suggest that we perform this computation step-by-step, through first-order infinitesimal extensions.  We will first set up this whole situation, starting with the short exact sequence of coefficients:
\[
0 \to \bigoplus_{\substack{\ell(J) = n \\ |J| = r}} K_* \to E_* / I^{r+1} \xrightarrow{C_r} E_* / I^r \to 0.
\]
Because we know that $E^\vee_* K(\Z/p^j, q)$ is a pro-free $E_*$-module deforming the original $K_* K(\Z/p^j, q)$, we tensor $E^\vee_* K(\Z/p^j, q)$ against the above short exact sequence to get a new short exact sequence:
\[
0 \to \bigoplus_J K_* K(\Z/p^j, q) \to E^\vee_* K(\Z/p^j, q) \otimes_{E_*} E_*/I^{r+1} \xrightarrow{H_r} E^\vee_* K(\Z/p^j, q) \otimes_{E_*} E_*/I^r \to 0.
\]
We can also build the free alternating coalgebraic ring $\Alt^* E^\vee_* B\Z/p^j$; tensoring the above short exact sequence of coefficients with any graded piece of this ring gives the exact sequence (which is not, a priori, left exact)
\[
\bigoplus_J \Alt^q K_* B\Z/p^j \to \Alt^q (E^\vee_* B\Z/p^j \otimes_{E_*} E_*/I^{r+1}) \xrightarrow{A_r} \Alt^q (E^\vee_* B\Z/p^j \otimes_{E_*} E_*/I^r) \to 0.
\]

Then, the cup product map $(B\Z/p^j)^{\wedge q} \to K(\Z/p^j, q)$ induces a map on homology \[\Alt^q E^\vee_* B\Z/p^j \xrightarrow{\circ} E^\vee_* K(\Z/p^j, q).\]  Bifunctoriality of the tensor product gives a diagram:
\begin{center}
\begin{tikzpicture}[
        dash line/.style={densely dotted},
        under line/.style={-stealth},
        normal line/.style={under line,
           preaction={draw=white, -, 
           line width=6pt}},
        onto/.style={normal line,->>},
        into/.style={normal line,>->},
        equal/.style={double,-},
    ]
    \matrix (m) [matrix of math nodes,
         row sep=2em, column sep=1.9em,
         text height=1.5ex,
         text depth=0.25ex]{
  & \bigoplus_J \Alt^q K_* B\Z/p^j & \Alt^q (E^\vee_* B\Z/p^j \otimes_{E_*} E_*/I^{r+1}) & \Alt^q (E^\vee_* B\Z/p^j \otimes_{E_*} E_*/I^r) & 0 \\
0 & \bigoplus_J K_* K(\Z/p^j, q) & E^\vee_* K(\Z/p^j, q) \otimes_{E_*} E_*/I^{r+1} & E^\vee_* K(\Z/p^j, q) \otimes_{E_*} E_*/I^r & 0. \\
    };
    \path[normal line]
        (m-1-2) edge node[right]{$\bigoplus_J (\circ \otimes \id_{K_*})$} (m-2-2)
                edge (m-1-3)
        (m-1-3) edge node[right]{$\circ \otimes \id_{E_*/I^{r+1}}$} (m-2-3)
                edge node[above]{$A_r$} (m-1-4)
        (m-1-4) edge node[right]{$\circ \otimes \id_{E_*/I^r}$} (m-2-4)
                edge (m-1-5)
        (m-2-1) edge (m-2-2)
        (m-2-2) edge node[above]{$I_r$} (m-2-3)
        (m-2-3) edge node[above]{$H_r$} (m-2-4)
        (m-2-4) edge (m-2-5);
\end{tikzpicture}
\end{center}

\begin{thm}
There is an isomorphism of graded coalgebraic rings \[\bigoplus_{q=0}^\infty \Alt^q E^\vee_* K(\Z/p^j, q) \cong \bigoplus_{q=0}^\infty E^\vee_* K(\Z/p^j, q).\]
\end{thm}
\begin{proof}
We perform an induction.  When $r = 1$, Ravenel and Wilson show that 
\[
\Alt^q (E^\vee_* B\Z/p^j \otimes_{E_*} E_* / I^1) = \Alt^q K_* B\Z/p^j \xrightarrow{\cong} K_* K(\Z/p^j, q) = E^\vee_* K(\Z/p^j, q) \otimes_{E_*} E_* / I^1
\]
is an isomorphism; this handles the base case.  This also tells us that the left-hand vertical map in the above diagram is an isomorphism.  In particular, both this map and the inclusion $I_r$ are injective, so their composite is injective, and hence so is the map $\bigoplus_J \Alt^q K_* B\Z/p^j \to \Alt^q (E/I^{r+1})_* B\Z/p^j$.  This means that the top sequence is short exact.

Then, assume that $\circ \otimes_{E_*} E_* / I^r$ induces an isomorphism for some fixed $r$, i.e., that the right-hand vertical map in the above diagram is an isomorphism of modules.  As the left-hand and right-hand vertical maps are isomorphisms, the center map must be as well.  As $q$ varies, the center maps additionally assemble into a map of graded coalgebraic rings, and so furthermore induce an isomorphism of rings.  So, induction provides ring isomorphisms for all $r$, and the Milnor sequence finishes the argument: \[E^\vee_* K(\Z/p^j, q) = \lim_r E^\vee_* K(\Z/p^j, q) \otimes_{E_*} E_*/I^r = \lim_r \Alt^q (E^\vee_* B\Z/p^j \otimes_{E_*} E_*/I^r) = \Alt^q E^\vee_* B\Z/p^j. \qedhere\]
\end{proof}

\section{Coalgebraic formal schemes, Cartier duals, and cohomology}

There are two relations between the formal schemes associated to $E$-cohomology and $E^\vee$-homology:
\begin{enumerate}
\item Let $A$ be an $E_*$-coalgebra whose underlying $E_*$-module is pro-free on a basis $\{e_i\}$.  If $I$ is a collection of indices such that $U_I = \Span \{e_i\}_{i \in I}$ is a subcoalgebra of $A$, $U_I$ is said to be a standard subcoalgebra.  Finally, the basis $\{e_i\}$ is said to be ``good'' if every finitely generated submodule of $A$ is contained in a standard subcoalgebra.  Coalgebras of this type have a functor to formal schemes given by the formula $\sch(A) = \colim_I \Spec U_I^\vee$, where $I$ runs over the indices corresponding to standard subcoalgebras.  The functor $\sch$ is an equivalence onto its essential image, with inverse given by $cX = \mathcal{O}_X^\vee$, and a formal scheme in its image is said to be a coalgebraic formal scheme.  Combining this language with the Hovey-Strickland results about $E$-theory, we see that \[\sch(E^\vee_* X) = \Spf E^* X\] for spaces $X$ with $K_* X$ even-concentrated.
\item Linear algebraic duality can also be encoded in Cartier duality.  For a formal group $G$ over a formal scheme $S = \colim S'$, we define its Cartier dual $DG$ by \[DG(R) := \underline{\Hom}_S(G, \G_m)(R) := \colim_{S'} \left\{(u, f) \middle| \begin{array}{c} u: \Spec R \to \Spec S', \\ f: u^* (G \times S') \to u^* (\G_m \times S')\end{array}\right\}.\]  For $X$ an even-concentrated $H$-space, so that $\Spf E^* X$ is a coalgebraic group, we have an isomorphism of $(\Spf E^*)$-schemes \[\Spf E^\vee_* X \cong D(\Spf E^* X).\]
\end{enumerate}

So far, we have seen that there is a function $(B\Z/p^\infty)_E^{\times q} \to K(\Z/p^\infty, q)_E$, which, if the scheme $(B\Z/p^\infty)_E^{\wedge q}$ exists, induces a homomorphism of groups $(B\Z/p^\infty)_E^{\wedge q} \to K(\Z/p^\infty, q)_E$.  Taking Cartier duals, we have a map \[\Spf \Alt^q E^\vee_* B\Z/p^\infty = D((B\Z/p^\infty_E)^{\wedge q}) \leftarrow DK(\Z/p^\infty, q)_E = \Spf E^\vee_* K(\Z/p^\infty, q),\] which in the previous section we showed to be an isomorphism.  We want to show now that the formal group scheme $(B\Z/p^\infty_E)^{\wedge q}$ does indeed exist and that the map $(B\Z/p^\infty_E)^{\wedge q} \to K(\Z/p^\infty, q)_E$ inducing the isomorphism on Cartier duals is itself an isomorphism.

Showing this amounts to showing that the coalgebras $\Alt^q E^\vee_* B\Z/p^\infty$ and $E^\vee_* K(\Z/p^\infty, q)$ admit good bases.  Every formal variety is a coalgebraic formal scheme, and it is further known that every first-order infinitesimal deformation of a formal variety to a formal scheme whose coordinate ring is free is again a formal variety.  Hence, we see that the coalgebras $E^\vee_* B\Z/p^\infty$ and $E^\vee_* K(\Z/p^\infty, q)$ have good bases.  We also have the following proposition of Strickland:
\begin{lem}(Strickland~\cite[Proposition 4.64]{StricklandFSFG})
Let $U_i$ be a diagram of coalgebras with colimit $U$, and assume that $U_i$ and $U$ all have good bases.  Then we have an interchange formula \[\colim \sch(U_i) = \sch(\colim U_i) = \sch U.\]
\end{lem}
Because the alternating coalgebraic group can be described as a colimit\footnote{This is a mess, but there is an $\mathbb{N}$-indexed diagram $F(cG) \times F(cG) \times F(cG^{q-1}) \times F(cG^q)^{n-1} \to F(cG^q)$ with four arrows:
\begin{enumerate}
\item Use $\ast$-multiplication on the two left-most factors in cG, then concatenate them with the third factor, then all the remaining factors in $F(cG^q)$.
\item Use the diagonal on the third factor, permute the third and second factors, concatenate the third and fourth factors, concatenate the first and second factors, then add everything in $F(cG^q)$.
\item Do (2), then apply $\chi$ and permute the first two factors of $cG$.
\item Do (2), then apply $\chi^q$ and a $q$-cycle permutation.
\end{enumerate}  Thinking of all these arrows as sharing a common target, the resulting colimit will encode the $q$th exterior power of the coalgebraic group $cG$.} and because $E^\vee_* B\Z/p^\infty$ is dual to a formal variety, we see that $\Alt^q E^\vee_* B\Z/p^\infty$ also has a good basis.  Hence the scheme $(B\Z/p^\infty_E)^{\wedge q}$ exists with the desired universal property, and the isomorphism from the previous section induces an isomorphism \[(B\Z/p^\infty_E)^{\wedge q} = \sch(\Alt^q E^\vee_* B\Z/p^\infty) \to \sch(E^\vee_* K(\Z/p^\infty, q) = K(\Z/p^\infty, q)_E.\]

\section{Conclusion}

There has been very recent buzz about Morava $E$-theory, concerning its construction as an $\mathcal{E}_\infty$-ring spectrum.  This problem was originally successfully tackled by Goerss, Hopkins, and Miller by the construction and analysis of an obstruction framework.  More recently, Lurie has announced the following result, of which the structured construction of Morava $E$-theory is the special case over a point:
\begin{thm}[Lurie, unpublished] (See Behrens-Lawson~\cite[Theorem 8.1.4]{BL})
Let $A$ be a local ring with $A / \mathfrak{m}_A$ a perfect field of positive characteristic $p$.  Then, let $X$ be a locally Noetherian, separated Deligne-Mumford stack over $\Spec A$ and $\mathbb{G}$ a $1$-dimensional $p$-divisible group over $X$ of constant height, and assume the existence of an \'etale cover $\pi: \tilde X \to X$ by a scheme $\tilde X$ such that for all points $x \in \tilde X^{\wedge}_{\mathfrak{m}_A}$ the map $\tilde X^\wedge_x \to \operatorname{Def}_{\pi^* \mathbb{G}_x}$ classifying the universal deformation of $(\pi^* \mathbb{G})|_{\tilde X^\wedge_x}$ is an isomorphism.

Then, given all this, there exists a sheaf of $\mathcal{E}_\infty$-ring spectra $E_{\mathbb{G}}$ on the \'etale site of $X^\wedge_{\mathfrak{m}_A}$, contravariantly assigned to the pair $(X, \mathbb{G})$, satisfying\ldots
\begin{itemize}
\item \ldots homotopy descent.
\item \ldots for any \'etale formal affine open $f: \Spf R \to X^\wedge_{\mathfrak{m}_A}$, the spectrum of sections $E_{\mathbb{G}}(R)$ is weakly even-periodic with coefficient ring $\pi_0 E_{\mathbb{G}}(R) = R$ and formal group described by an isomorphism $\gamma_f: f^* \mathbb{G}^0 \to \CP^\infty_{E_{\mathbb{G}}(R)}$ naturally in $f$. \qed
\end{itemize}
\end{thm}

A next goal, then, is to compute the $E_{\mathbb{G}}(R)$-homology of Eilenberg-Mac Lane spaces, which is likely now approachable by using these local computations in Morava $E$-theory and investigating the homotopy descent property of the sheaf of spectra $E_{\mathbb{G}}$.

\newpage

\bibliography{e-thy}
\bibliographystyle{plain}

\end{document}